\documentstyle[leqno]{bcp01e}
\def\operatorname#1{{\rm#1\,}}
\def\text#1{{\hbox{#1}}}
\def\range{\operatorname{Range}}
\def\span{\operatorname{Span}}
\def\rank{\operatorname{Rank}}
\def\id{\operatorname{Id}}
\def\Gr{\operatorname{Gr}}
\newcommand{\reals}{\mbox{${\rm I\!R }$}}
\begin{document}
\mathclass{Primary 53B30; Secondary 53C05.}

\newtheorem{th}{Theorem}
\newtheorem{lemma}[th]{Lemma}

\abbrevauthors{P. Gilkey and R. Ivanova}
\abbrevtitle{Complex IP curvature tensors}

\title{Complex IP pseudo-Riemannian\\
algebraic curvature tensors}

\author{Peter\ B\ Gilkey}
\address{Mathematics Department, University of Oregon,
Eugene Or 97403 USA\\
E-mail: gilkey@darkwing.uoregon.edu}

\author{Raina\ Ivanova}
\address{Dept. of Descriptive Geometry, University of Architecture,\\
Civil Engineering \& Geodesy, 1, Christo Smirnenski Blvd., 1421 Sofia,
Bulgaria\\email:  ivanovar@is.tsukuba.ac.jp}

\maketitlebcp

\def\operatorname#1{{\rm#1\,}}
\def\text#1{{\hbox{#1}}}

\section{Introduction}\label{Sect1} The Riemann curvature tensor contains a great
deal of information about the geometry of the underlying pseudo-Riemannian
manifold; pseudo-Riemannian geometry is to a large extent the study of this
tensor and its covariant derivatives. It is often convenient to work in a purely
algebraic setting. We shall say that a tensor is an algebraic curvature tensor if
it satisfies the symmetries of the Riemann curvature tensor. The Riemann
curvature tensor defines an algebraic curvature tensor at each point of the
manifold; conversely every algebraic curvature tensor is locally geometrically
realizable. Thus algebraic curvature tensors are an integral part of certain questions in
differential geometry.

The skew-symmetric curvature operator is a natural object of study; there are
other natural operators and we refer to \cite{PI5} for a survey of this area. In
the present paper, we examine when the (complex) Jordan normal form of the
skew-symmetric curvature operator is constant either in the real or complex
settings. In Section
\ref{Sect2}, we give the basic definitions and notational conventions we shall
need. We also review the basic results in the real setting. In Section
\ref{Sect3}, we discuss a natural generalization to the complex setting.

\section{IP algebraic curvature tensors}\label{Sect2}
Let $(M,g)$ be a connected pseudo-Riemannian manifold of signature $(p,q)$ and
dimension $m=p+q$. Let $\nabla$ be the Levi-Civita connection. The curvature
operator and associated curvature tensor are defined by:
\begin{eqnarray}\label{2.1}
&&{}^gR(X,Y):=\nabla_X\nabla_Y-\nabla_Y\nabla_X-\nabla_{[X,Y]},\\
&&{}^gR(X,Y,Z,W):=g({}^gR(X,Y)Z,W).\nonumber\end{eqnarray}
This tensor has the following symmetries:
\begin{eqnarray}  \label{eqna}
     &&{}^gR(X,Y,Z,W)=-{}^gR(Y,X,Z,W),\\
     &&{}^gR(X,Y,Z,W)={}^gR(Z,W,X,Y),\text{ and}\label{eqnb}   \\
     &&{}^gR(X,Y,Z,W)+{}^gR(Y,Z,X,W)+{}^gR(Z,X,Y,W)=0.\label{eqnc}
\end{eqnarray}

Let $V$ be finite dimensional real vector space which is
equipped with a non-degenerate inner product $(\cdot,\cdot)$ of signature $(p,q)$.
We say that a $4$ tensor $R\in\otimes^4V$ is an  {\it algebraic curvature tensor}
if
$R$ satisfies the symmetries given in equations (\ref{eqna}-\ref{eqnc}). 
Note that we do not impose the second Bianchi identity; algebraic curvature
tensors measure second order phenomena. We say that $(M,g)$ is a
{\it geometric realization} of an algebraic curvature tensor at a point $P$ if
there is an isometry $\Psi:T_PM\rightarrow V$ so that ${}^gR_P=\Psi^*R$.

Let $\{e_1,e_2\}$ be an oriented basis for a $2$ plane $\pi\subset V$ and let
$h_{ij}:=(e_i,e_j)$ describe the restriction of the inner product $(\cdot,\cdot)$
to
$\pi$. We shall say that
$\pi$ is {\it non-degenerate} if $h_{ij}$ is non-degenerate, i.e.
$\det(h):=h_{11}h_{22}-h_{12}^2\ne0$. We say that $\pi$ is {\it timelike}, {\it mixed}, or
{\it spacelike} as the quadratic form
$h_{ij}$ is negative definite, indefinite, or positive definite, respectively. We
may decompose the Grassmannian of all oriented non-degenerate $2$ planes
as the disjoint union of the oriented timelike, mixed, and spacelike $2$ planes.

Let $R$ be an algebraic curvature tensor on $V$. Let $\{e_1,e_2\}$ be an oriented
basis for a non-degenerate spacelike $2$ plane $\pi\subset V$. The {\it skew-symmetric
curvature operator}
$$R(\pi):=\det(h)^{-1/2}R(e_i,e_j)$$
is independent of the particular oriented basis for $\pi$ which was chosen. We say
that $R$ is {\it spacelike Jordan IP} if the complex Jordan normal form of the operator
$R(\pi)$ is constant on the Grassmannian  $\Gr_{0,2}^+(V)$ of oriented spacelike $2$ planes. The notions
of {\it timelike Jordan IP} and {\it mixed Jordan IP} are defined similarly using the
Grassmannians $\Gr_{2,0}^+(V)$ and $\Gr_{1,1}^+(V)$ respectively. In the
Riemannian setting $(p=0)$, this is equivalent to assuming that the eigenvalues
of $R(\pi)$ are constant on $\Gr_{(0,2)}^+(V)$. In the pseudo-Riemannian setting ($p>0$) a
bit more care must be taken as there are examples where $R(\pi)$ has only the
zero eigenvalue but where the rank of $R(\pi)$ varies with $\pi$; we refer to
\cite{GZ6} for details. We note that there are examples of algebraic
curvature tensors which are spacelike Jordan IP but not timelike Jordan IP; again see
\cite{GZ6}. We say that a pseudo-Riemannian manifold
$(M,g)$ is {\it spacelike Jordan IP} if the associated curvature tensor ${}^gR$ is
spacelike Jordan IP at every point of $M$; the eigenvalues and Jordan normal form are
allowed to vary with the point. 

Let $p=0$. The study of the skew-symmetric curvature
tensor was initiated in this context by Stanilov and Ivanova
\cite{IS8}; we also refer to related work by Ivanova
\cite{I9}-\cite{I13}. Subsequently, Ivanov and Petrova \cite{IP7} classified
the spacelike Jordan IP metrics in the Riemannian setting for $m=4$; for this reason the notation `IP'
has been used by later authors. The classification in \cite{IP7} was later extended by P. Gilkey, J.
Leahy, and H. Sadofsky
\cite{GLS3} and by Gilkey
\cite{G1} to the cases $m\ge5$ and
$m\ne7$. We refer to Gilkey and Semmelman \cite{GS4} for some partial results if
$m=7$. Zhang \cite{Z12} has extended these results to the Lorentzian setting
($p=1$); see also Gilkey and Zhang \cite{GZ7} for related work.

We say that $R$ has {\it spacelike rank $r$} if $\rank(R(\pi))=r$ for every
spacelike $2$ plane $\pi$. The following theorem, which shows that $r=2$ in many cases, was proved using
topological methods \cite{GLS3,Z12}.
\begin{th}\label{Th1}
Let
$R$ be an algebraic curvature tensor of spacelike rank $r$ on a vector
space of signature $(p,q)$.
\smallbreak{\rm 1)} Let $p\le1$. Let $q=5$, $q=6$, or
$q\ge9$. Then $r=2$.
\smallbreak{\rm 2)} Let $p=2$. Let $q\ge10$. Assume neither $q$ nor $q+2$ are
powers of $2$. Then
$r=2$.
\end{th}

In light of Theorem \ref{Th1}, we shall focus our attention on the algebraic
curvature tensors of rank $2$. Let $V$ be a vector space of signature $(p,q)$.
We say that a linear map $\phi$ of $V$ is {\it admissible} if $\phi$ satisfies the following
two conditions:
\begin{enumerate}
\item $\phi$ is self-adjoint and $\phi^2=\id$, or $\phi^2=-\id$, or $\phi^2=0$.
\item If $\phi(x)=0$, then $(x,x)\le0$, i.e. $\ker\phi$ contains no spacelike vectors.
\end{enumerate}

Let $\phi$ be admissible. If $\phi^2=\id$, then $\phi$ is an isometry, i.e. $(\phi
x,\phi y)=(x,y)$ for all
$x,y\in V$. If $\phi^2=-\id$, then $\phi$ is a {\it para-isometry}, i.e. $(\phi x,\phi y)=(-x,-y)$ for
all $x,y\in V$; necessarily $p=q$ in this setting. If $\phi^2=0$, then the range of $\phi$ is
{\it totally isotropic}, i.e. $(\phi x,\phi y)=0$ for all $x,y\in V$. We define:
\begin{eqnarray*}
   &&R_\phi(x,y)z:=(\phi y,z)\phi x
   -(\phi x,z)\phi y,\text{ and}\label{RP}\\
   &&R_\phi(x,y,z,w):=(\phi y,z)(\phi x,w)
   -(\phi x,z)(\phi y,w).\nonumber\end{eqnarray*}
We showed in \cite{GZ6} that $R_\phi$ is an algebraic curvature tensor with
$range(R_\phi(\pi))=\phi\pi$ if $\pi$ is a spacelike $2$ plane. The eigenvalue structure of
$R_\phi(\pi)$ is given by:
\begin{enumerate}
\item Suppose that $\phi^2=\pm\id$. Then $R_\phi(\pi)$ is a rotation through an
angle of 90 degrees on the spacelike ($\phi^2=\id$) or timelike ($\phi^2=-\id$) $2$ plane
$\phi\pi$,
$R_\phi(\pi)$ vanishes on $\phi\pi^\perp$, and $R_\phi$ has two
non-trivial complex eigenvalues $\pm\sqrt{-1}$.
\item Suppose that $\phi^2=0$. Since $\ker\phi$ contains no spacelike vectors,
$R_\phi(\pi)$ has rank $2$. The $2$ plane $\phi\pi$ is totally isotropic. We have
$\{R_\phi(\pi)\}^2=0$.
\end{enumerate}
Thus $CR_\phi$ is a spacelike rank $2$ Jordan IP algebraic curvature tensor for any $C\ne0$. Conversely,
we have the following classification result \cite{GZ6}:
\begin{th}\label{Thm2} Let $V$ be a vector space of signature $(p,q)$, where $q\ge5$. A tensor $R$ is
a spacelike rank 2 Jordan IP algebraic curvature tensor on $V$ if and only if there exists a non-zero
constant
$C$ and an admissible $\phi$ so that
$R=CR_\phi$.
\end{th}

\section{ Almost complex Jordan IP algebraic curvature tensors}\label{Sect3}
Algebraic curvature tensors have been studied by many authors in the complex
setting; we refer to Falcitelli, Farinola, and Salamon
\cite{FFS} and Gray \cite{Gr} for further details concerning almost
Hermitian geometry. 

Let $J:V\rightarrow V$ be a real linear map with $J^2=-\id$. We use $J$ to provide
$V$ with a complex structure:
$(a+\sqrt{-1}b)v=av+bJv$. Thus a real linear map $S$ of $V$ is {\it complex} if
and only if  $SJ=JS$. We shall assume that $J$ is {\it pseudo-Hermitian}, i.e.
$(Jv_1,Jv_2)=(v_1,v_2)$; necessarily both $p$ and $q$ are even.

A  $2$ plane $\pi$ is called  {\it complex line} if and only if
$J\pi\subset\pi$. If $\pi$ is non-degenerate, then $\pi$ is either spacelike or
timelike, there are no mixed complex lines. 

An algebraic curvature tensor $R$ is said to be {\it
almost complex} if $JR(x,Jx)=R(x,Jx)J$ for all $x$ in $V$, i.e. $R(x,Jx)$ is complex linear.
Such an
$R$ is said to be {\it almost complex spacelike Jordan IP} if $R(\pi)$ (regarded as a complex linear
map) has constant Jordan normal form for every spacelike complex line; the notion of {\it almost complex
timelike Jordan IP} is defined similarly. 

Theorem \ref{Thm2} controls the eigenvalue
structure of a spacelike Jordan IP algebraic curvature tensor in the real setting. There
is a similar result in the complex setting which we describe as follows.
Let $p=0$ and let $R$ be an almost complex
spacelike Jordan IP algebraic curvature tensor. The operator
$JR(\pi)$ is a self-adjoint complex linear map and is therefore diagonalizable.
Let
$\{\lambda_i,\mu_i\}$ be the eigenvalues and multiplicities of $JR(\cdot)$,
where $\mu_0\ge\mu_1\ge...\ge\mu_\ell$; $\lambda_i\in\reals$  and $q=2(\mu_0+...+\mu_\ell)$. We refer to
\cite{G2} for the proof of the following result which controls the eigenvalue
structure in the Riemannian setting; it is not known if a similar result holds in
the higher signature setting.
\begin{th}\label{thm3} Let $R$ be an almost complex spacelike Jordan IP algebraic
curvature tensor on a Riemannian vector space of signature $(0,q)$. Let
$\{\lambda_i,\mu_i\}$ be the eigenvalues and multiplicities of $JR(\cdot)$,
where $\mu_0\ge...\ge\mu_\ell>0$. Suppose
$\ell\ge1$. If
$m\equiv2$ mod $4$, then $\ell=1$ and $\mu_1=1$. If $m\equiv0$ mod $4$, then
either
$\ell=1$ and $\mu_1\le2$ or
$\ell=2$ and $\mu_1=\mu_2=1$.
\end{th}

We say that $(\phi,J)$ is an {\it admissible pair} if
$\phi$ is admissible, if $\phi J=\pm J\phi$, and if $J$ is a pseudo-Hermitian almost complex structure
on $V$.

\begin{th}\label{Thm4} Let $V$ be a vector space of signature $(p,q)$. If $(\phi,J)$ is an admissible
pair, then
$R_\phi$ is an almost complex spacelike Jordan IP algebraic curvature tensor.
\end{th}

\Proof By Theorem \ref{Thm2}, $R_\phi$ is a spacelike rank $2$ Jordan IP algebraic curvature tensor.
Since $J\phi=\varepsilon\phi J$ for $\varepsilon=\pm1$, we may compute:
\begin{eqnarray*}
&&JR_\phi(x,Jx)z=(\phi Jx,z)J\phi x-(\phi x,z)J\phi J x,\\
&&R_\phi(x,Jx)Jz=(\phi Jx,Jz)\phi x-(\phi x,Jz)\phi Jx\\
&&\quad=-(J\phi Jx,z)\phi x+(J\phi x,z)\phi Jx\\
&&\quad=\varepsilon^2(\phi Jx,z)J\phi x-\varepsilon^2(\phi x,z)J\phi Jx
   =JR_\phi(x,Jx)z.
\end{eqnarray*}
This shows that $R$ is almost complex. \endproof

Let $V$ be a vector space of signature $(p,q)$. We say that $(\phi_1,\phi_2,J)$ is an {\it admissible
triple} if the following conditions are satisfied:
\begin{enumerate}
\item $\phi_1$ and $\phi_2$ are admissible and either $\phi_1^2\ne0$ or $\phi_2^2\ne0$;
\item $J$ is a pseudo-Hermitian almost complex structure on $V$;
\item $J\phi_1=\phi_1J$, $J\phi_2=-\phi_2J$, and $\phi_2\phi_1+\phi_1\phi_2=0$.
\end{enumerate}

\begin{lemma}\label{Lem3} Let $V$ be a vector space of signature $(p,q)$. Let $\{\phi_1,\phi_2,J\}$ be
an admissible triple on $V$.\begin{enumerate}
\item If $x$ is spacelike, then the set $\{\phi_1x,\phi_1Jx,\phi_2x,\phi_2Jx\}$ is
orthogonal and linearly independent.
\item For any $x\in V$, we have that:
$$R_{\phi_1}(x,Jx)R_{\phi_2}(x,Jx)=R_{\phi_2}(x,Jx)R_{\phi_1}(x,Jx)=0.$$
\end{enumerate}\end{lemma}

\Proof To show that $\{\phi_1x,\phi_1Jx,\phi_2x,\phi_2Jx\}$ is an orthogonal set, we compute:
\begin{eqnarray*}
&&(\phi_1x,\phi_1Jx)=\varepsilon_1(\phi_1x,J\phi_1x)
    =-\varepsilon_1(J\phi_1x,\phi_1x)=-(\phi_1Jx,\phi_1x),\\
&&(\phi_1x,\phi_2x)=(\phi_2\phi_1x,x)=-(\phi_1\phi_2x,x)=-(\phi_2x,\phi_1x),\\
&&(\phi_1x,\phi_2Jx)=-(J\phi_2\phi_1x,x)
     =\varepsilon_1\varepsilon_2(\phi_1\phi_2Jx,x)=-(\phi_2Jx,\phi_1x),\\
&&(\phi_1Jx,\phi_2x)=-(x,J\phi_1\phi_2x)
      =\varepsilon_1\varepsilon_2(x,\phi_2\phi_1Jx)=-(\phi_2x,\phi_1Jx),\\
&&(\phi_1Jx,\phi_2Jx)=(\phi_2\phi_1Jx,Jx)=-(\phi_1\phi_2Jx,Jx)
     =-(\phi_2Jx,\phi_1Jx),\\
&&(\phi_2x,\phi_2Jx)=\varepsilon_2(\phi_2x,J\phi_2x)
    =-\varepsilon_2(J\phi_2x,\phi_2x)=-(\phi_2Jx,\phi_2x).
\end{eqnarray*}

Let $\pi:=\span\{x,Jx\}$. Then $\phi_1\pi$ and $\phi_2\pi$ are orthogonal $2$ planes. If
$\phi_1^2=\pm\id$, then $\phi_1\pi$ is non-degenerate so $\phi_2\pi\subset\phi_1\pi^\perp$
implies $\phi_1\pi\cap\phi_2\pi=\{0\}$ and assertion (1) follows; the argument is the same if
$\phi_2^2=\pm\id$. 

To prove assertion (2), we compute:
\begin{eqnarray*}
  &&R_{\phi_1}(\pi)R_{\phi_2}(\pi)z
  =R_{\phi_1}(\pi)\{(\phi_2Jx,z)\phi_2x-(\phi_2x,z)\phi_2Jx\}\\
  &&\quad=(\phi_2Jx,z)\{(\phi_1Jx,\phi_2x)\phi_1x-(\phi_1x,\phi_2x)\phi_1Jx\}\\
  &&\quad-(\phi_2x,z)\{(\phi_1Jx,\phi_2Jx)\phi_1x-(\phi_1x,\phi_2Jx)\phi_1Jx\}
     =0.
\end{eqnarray*}
We argue similarly to show that $R_{\phi_2}(\pi)R_{\phi_1}(\pi)=0$.
\endproof

The following is the main result of this paper. It shows the estimates
of Theorem
\ref{thm3} are sharp and provides a large family of non-trivial new examples.

\begin{th}\label{Thm5} Let $V$ be a vector space of signature $(p,q)$. Let $(\phi_1,\phi_2,J)$ be an
admissible triple on $V$. Let
$\lambda_i$ be real constants. Then
$R:=\lambda_1R_{\phi_1}+\lambda_2R_{\phi_2}$ is an almost complex spacelike Jordan IP
algebraic curvature tensor.
\end{th}

\Proof We assume $\lambda_1\ne0$ and $\lambda_2\ne0$;
otherwise the proof follows directly from Theorem \ref{Thm4}. As the set of almost complex
algebraic curvature tensors is a linear subspace of the set of all $4$ tensors,
Theorem \ref{Thm4} shows that $R$ is an almost complex algebraic
curvature tensor. 

We complete the proof by discussing the complex Jordan form. Let
$\{x,Jx\}$ be an orthonormal basis for a spacelike complex line $\pi$. We use Lemma \ref{Lem3} to see
that $\phi_1\pi$ and $\phi_2\pi$ are orthogonal complex lines and that $\rank(R(\pi))=4$. Also by Lemma
\ref{Lem3} we have:
$$R(\pi)^2=\lambda_1^2R_{\phi_1}(\pi)^2+\lambda_2^2R_{\phi_2}(\pi)^2.$$

 Suppose that $\phi_1^2=\pm\id$ and that $\phi_2^2=\pm\id$. The
metric restricted to
$\phi_1\pi\oplus\phi_2\pi$ is non-degenerate. Let
$$V_0:=(\phi_1\pi\oplus\phi_2\pi)^\perp.$$We have an orthogonal direct sum
decomposition
$$V=\phi_1\pi\oplus\phi_2\pi\oplus V_0$$ which is preserved by $R(\pi)$. The map
$R(\pi)$ has 4 non-trivial purely-imaginary eigenvalues $\pm\lambda_i\sqrt{-1}$;
$V_0=\ker(R(\pi))$. The map $R(\pi)^2$ is diagonalizable. Thus
$R$ is almost complex spacelike Jordan IP because: 
\begin{eqnarray*}
  &&R(\pi)^2=0\text{ on }V_0,\\
  &&R(\pi)^2=-\lambda_1^2\text{ on }\phi_1\pi\text{, and}\\
  &&R(\pi)^2=-\lambda_2^2\text{ on }
  \phi_2\pi.\end{eqnarray*}

 Suppose that $\phi_1^2=\pm\id$ and $\phi_2^2=0$; the argument is similar if
$\phi_1^2=0$ and $\phi_2^2=\pm\id$. By assumption
$\ker\phi_2$ contains no spacelike vectors. Let
$V_0:=(\phi_1\pi)^\perp$. Then we have an orthogonal direct sum decomposition $V=\phi_1\pi\oplus V_0$
which is preserved by $R(\pi)$. The map $R(\pi)$ has two non-zero eigenvalues
$\pm\lambda_1\sqrt{-1}$. The map $R(\pi)^2$ is diagonalizable; 
$$R(\pi)^2=0\text{ on }V_0\text{ and }R(\pi)^2=-\lambda_1^2\text{ on }\phi_1\pi.$$
Thus $R$ is almost
complex spacelike Jordan IP.\endproof

We construct examples to show that all 8 cases of
Theorem \ref{Thm5} can occur. Let
\smallbreak\noindent\begin{eqnarray*}
  &&e_1:=\left(\begin{array}{rrrr}0&1&0&0\\
      1&0&0&0\\0&0&0&1\\0&0&1&0
      \end{array}\right)\text{ and } 
  e_2:=\left(\begin{array}{rrrr}0&0&1&0\\
      0&0&0&-1\\1&0&0&0\\0&-1&0&0\end{array}\right)\text{ on }\reals^{(0,4)},\\ \\
 &&J_0:=\left(\begin{array}{rr}0&1\\-1&0\end{array}\right)
   \qquad\phantom{..}\text{ and }
  \alpha:=\left(\begin{array}{rr}0&1\\1&0\end{array}\right)
   \qquad\qquad\phantom{....}\text{ on }
   \reals^{(0,2)},\\ \\
 &&\beta:\phantom{}=\left(\begin{array}{rr}0&1\\-1&0\end{array}\right)
\qquad\phantom{..}
\text{ and }\gamma:=\left(\begin{array}{rr}1&-1\\1&-1\end{array}\right)
\qquad\qquad\phantom{.}
\text{ on }\reals^{(1,1)}.\end{eqnarray*}
be matrices satisfying the relations\medbreak
\centerline{\begin{tabular}{rrrrr}
$e_1^*=e_1,$&$e_2^*=e_2,$&$e_1^2=\id,$&$e_2^2=\id,$&$e_1e_2+e_2e_1=0$,\\ \\
$J_0^*=-J_0,$&$J_0^2=-\id,$&$\alpha^*=\alpha,$&$\alpha^2=\id,$&$J_0\alpha=-\alpha J_0,$\\ \\
$\beta^*=\beta,$&$\beta^2=-\id,$&$\gamma^*=\gamma,$&$\gamma^2=0$,&$\range\gamma=\ker\gamma.$
\end{tabular}}\medbreak\noindent

Define the matrix $\tau_i$ by:
$$\tau_i=\left\{\begin{array}{ll}
\id\otimes\id&\quad\text{ if }\delta_i=+1,\\
\beta\otimes\id&\quad\text{ if }\delta_i=-1,\\
\id\otimes\gamma&\quad\text{ if }\delta_i=0.\end{array}\right.$$

We construct
$(\phi_1,\phi_2,J)$ admissible so 
$$\phi_1^2=\delta_1\id\text{ and }
\phi_2^2=\delta_2\id$$ by setting:
\begin{eqnarray*}
  &&\phi_1:=e_1\otimes\id\otimes\tau_1,\\ 
  &&\phi_2:=e_2\otimes\alpha\otimes\tau_2,\text{ and }\\
  &&J:=\id\otimes J_0\otimes\id.\end{eqnarray*}
The tensor $$R=\lambda_1R_{\phi_1}+\lambda_2R_{\phi_2}$$ is then both
almost complex spacelike Jordan IP and
almost complex timelike Jordan IP.

\section{Acknowledgements} The research of the first author was partially
supported by the NSF (USA). The research of the second author was partially
supported by the JSPS Post Doctoral Fellowship Program (Japan). The research of
both authors was also partially supported and facilitated by a joint
visit to the MPI (Leipzig, Germany).

\end{document}